\documentclass{amsart}

\usepackage{amssymb}

\newtheorem{thm}{Theorem}[section] 
\newtheorem{lm}[thm]{Lemma}
\newtheorem{cor}[thm]{Corollary}

\newtheorem{prop}[thm]{Proposition}
\DeclareMathOperator{\aut}{Aut}
\newcounter{szam}
\newcommand{\pf}{\noindent {\it Proof.} \;}

\begin{document}

\title{Group invariants of certain Burn loop classes}
\author{G\'abor P. Nagy}
\thanks{Supported by OTKA Grants Nos. F021271 and T020066.}
\thanks{Received in August 1997}
\thanks{Communicated by Frank De Clerck}
\address{JATE Bolyai Institute, Aradi v\'ertan\'uk tere 1. H-6720 Szeged
(Hungary)} 
\email{nagyg@math.u-szeged.hu}

\maketitle

\section{Introduction}
To any loop $(L,\cdot)$, one can associate several groups, for example its
multiplication groups $G_\mathrm{left}(L)$ and 
$G_\mathrm{right}(L)$ and $M(L)=\langle G_\mathrm{left}(L),
G_\mathrm{right}(L) \rangle$, the groups of (left or right)
pseudo-automorphisms, group of automorphisms, or the group of collineations
of the associated 3-net. Groups, which are isotope invariants are of
special interest. For example, the groups $G_\mathrm{left}(L)$, 
$G_\mathrm{right}(L)$ and $M(L)$ are isotope invariant for any loop $L$.
These groups contain many information about the loop $L$, 
the standard references on this field are \cite {BarlStr}, \cite{Bruck},
\cite{Pflug}. 

For some special loop classes, other isotope invariant groups can be
defined. For Bol loops, M. Funk and P.T. Nagy \cite{Funk-Nagy} investigated {\em
the collineation group generated by the Bol reflections}. The notion
of the {\em core} was first studied by R.H. Bruck \cite{Bruck} for Moufang
loops and by V.D. Belousov \cite{Belousov} for Bol loops. Recently, this
concept was intensively used in P.T. Nagy and K. Strambach \cite{NStisrael}. 

In the papers \cite{Burn1,Burn2}, R.P. Burn defined the infinite classes
$B_{4n}$, $n\leq 2$ and $C_{4n}$, $n\leq 2$, $n$ even of Bol loops. These
examples satisfy the left conjugacy closed property, that is, their section
\[S(L) = \{ \lambda_x : x \in L\}\] is
invariant under conjugation with elements of the group
$G_\mathrm{left}(L)= \langle \lambda_x | x \in L \rangle$
generated by the (left) translations $\lambda_x: y \mapsto xy$. 

In this paper, we determine the collineation groups generated by the Bol
reflections, the core, the automorphism groups and the full direction
preserving collinea\-tion groups of the loops $B_{4n}$ and $C_{4n}$ given by
R.P. Burn. We also prove some lemmas and use new methods in order to
simplify the calculations in these groups.
\section{Basic concepts}
A loop $L$ is said to be a {\em Bol loop}, if 
\[x\cdot (y \cdot xz) = (x \cdot yx) \cdot z\]
holds for all $x,y,z \in L$. This is equivalent with $\lambda_x \lambda_y
\lambda_x \in S(L)$ for all $x,y \in L$. In any Bol loop, the group
\begin{equation} \label{refl}
N = \langle (\lambda_x^{-1} \rho_x^{-1}, \lambda_x) | x \in L \rangle
\end{equation}
is a normal subgroup of the directions preserving collineation group of the
3-net belonging to the loop $L$, cf. \cite{Funk-Nagy}, \cite{Ngnote}.
Actually, the fact that $(\lambda_x^{-1} \rho_x^{-1}, \lambda_x)$ is a
direction preserving collineation for all $x \in L$ is equivalent with the
Bol property for the coordinate loop. As in \cite{Funk-Nagy}, we define the
endomorphism $\Phi$ by 
\begin{equation} \label{Ndef}
\Phi: \left \{
\begin{array}{l}
N \to G(L)\\
(\lambda_x^{-1} \rho_x^{-1}, \lambda_x) \mapsto \lambda_x.
\end{array} \right.
\end{equation}
This map $\Phi$ will help us to determine the group $N$, which acts
transitively on the set of horizontal lines, and so, plays an important
role in the description of the full collineation group of the 3-net. In
general, about the kernel of $\Phi$ one can only know, that it is
isomorphic to a subgroup of the left nucleus of $L$ (see \cite{Funk-Nagy},
Theorem 3.1).

The {\em core} of a Bol loop $(L,\cdot)$ is the groupoid $(L,+)$, where the
binary operation ``$+$'' is defined by 
\[x+y = x \cdot y^{-1} x, \hskip 1cm x,y \in L.\]
This groupoid satisfies the following identities:
\[ \begin{array}{ll}
x+x=x &\\
x+(x+y)=y 
&\forall x,y,z \in L \\
x+(y+z)=(x+y)+(x+z)&
\end{array} \]
An altarnative way to define the core is via the action of the Bol
reflections on the set of vertical lines of the associated 3-net. In this
way, the core turns out to be strongly related to the group $N$. 

We say that the loop $L$ is {\em left conjugacy closed}, if $S(L)$ is
invariant under the conjugation with the elements of $G(L)$. This concept
was introduced in the paper \cite{NStcanad} by P.T. Nagy and K. Strambach. 
They also defined the notion of {\em Burn loop}, which is a left conjugacy
closed Bol loop. Examples for such loops are the following constructions
due tu R.P. Burn \cite{Burn1,Burn2}. 

The section $S(L)$ of a loop $L$ is a sharply transitive set of
permutations, for any $x \in L$, there is a uniquely defined $\lambda_x$
mapping the unit element $1$ to $x$. Thus, by $x \cdot y= y^{\lambda_x}$,
the multiplication of $L$ is given by the set $S(L)$ and the
unit elements $1$. Theorem I.7 in \cite{Burn1,Burn2} says that if the set
$S(L)$ is invariant under conjucation with its own elements, different
choises of the unit element still give isomorphic loops, hence a Burn loop
is completely determined by its section $S(L)$ (up to isomorphism).

{\bf The loops $B_{4n}$, $n\geq 2$:} Let the group $G_{8n}$ be generated by
the
elements $\alpha$, $\beta$, $\gamma$ with the relations $\alpha^{2n} =
\beta^2 = \gamma^2 = (\alpha \beta)^2 = id$, $\alpha \gamma = \gamma
\alpha$ and $\beta \gamma = \gamma \beta$. 

The set $S(B_{4n})$ will be
\[ S(B_{4n}) = \{ \alpha^{2i}, \alpha^{2j+1} \beta, \alpha^k \beta \gamma :
i,j \in \{1, \ldots, n\}, k \in \{1, \ldots, 2n \}\}. \]
Then, with the action of $G_{8n}$ on the right cosets of $\langle \beta
\rangle$, $B_{4n}$ is a Burn loop, for all $n \geq 2$. Moreover, it is
non-Moufang, see \cite{Burn1,Burn2}.

{\bf The loops $C_{4n}$, $n\geq 2$, $n$ even:} Let the group $H_{8n}$ be
\[H_{8n} = \langle \alpha, \beta, \gamma : \alpha^{2n} = \beta^2 = \gamma^2
= (\alpha \beta)^2 = id, \alpha\gamma = \gamma \alpha, \beta \gamma =
\gamma \beta \alpha^n \rangle. \]

The set $S(C_{4n})$ will be
\[ S(C_{4n}) = \{ \alpha^{2i}, \alpha^{2j+1} \beta, \alpha^k \beta \gamma :
i,j \in \{1, \ldots, n\}, k \in \{1, \ldots, 2n \}\}. \]
Choosing again the action of $H_{8n}$ on the right cosets of $\langle \beta
\rangle$, the section $C_{4n}$ becomes a Burn loop, for all $n \geq 2$, $n$
even. It is non-Moufang, see \cite{Burn1,Burn2}.

In \cite{NStcanad}, the authors showed that the square of any
element of a Burn loop belongs to the intersection of the left
and middle nuclei. In any Bol loop, these two
nuclei coincide (cf. \cite{Ngnote}, Proposition 2.1) and form a
normal subgroup of the loop (see Lemma \ref{normnuc}). Thus, if 
$L$ denotes a (left) Bol loop, one can speak of the factor loop
$L/N_\lambda$. 
\begin{lm} \label{normnuc}
Let $(L,\cdot)$ be a (left) Bol loop. Then its left nucleus
$N_\lambda$ is a normal subgroup of $L$. 
\end{lm}
\pf Let $(L,\cdot)$ be a left Bol loop. Let us denote by 
$G_{\mathrm{left}}(L)$ and $G_{\mathrm{right}}(L)$ the groups 
generated by the left and right translations of $L$, respectively. 
Let $M(L)$ denote the group generated by $G_{\mathrm{left}}(L)$ 
and $G_{\mathrm{right}}(L)$. The Bol identity $x\cdot (y \cdot xz) 
= (x\cdot yx) \cdot z$ can also be expressed by $\rho_{xz} \lambda_x
= \rho_x \lambda_x \rho_z$, or equivalently, $\lambda_x \rho_z 
\lambda_x^{-1} = \rho_{xz} \rho_z^{-1} \in G_{\mathrm{right}}(L)$. 
This means that $G_{\mathrm{right}}(L)$ is a normal subgroup of $M(L)$.

Let now $u$ be a permutation of $L$ with $1^u=n$ and let us suppose that
$u$ centralizes the group $G_{\mathrm{right}}(L)$. Then we have for any 
$x \in L$
\[x^u = 1^{\rho_x u} = 1^{u \rho_x} = nx,\]
that is, $u=\lambda_n$. Moreover, $\lambda_n \rho_x = \rho_x \lambda_n$
for all $x \in L$
means exactly that $n$ is an element of the left nucleus $N_\lambda(L)$
of $L$. 
Hence, $U=\{ \lambda_n: n\in N_\lambda(L)\}$ is the centralizer of
the normal subgroup $G_{\mathrm{right}}(L)$ in $M(L)$, it should also be 
normal. This implies that $N_\lambda(L) = 1^U$ is a normal subgroup
of $L$, see \cite{Albert1}, Theorem 3.\qed
{\it Remark.} Clearly, if $L$ is a Burn loop, the factor loop
$L/N_\lambda$ is Burn as well. This means that in the quotient
loop $L/N_\lambda$ of a Burn loop $L$ every element has order 2.
\section{The kernel of the map $\Phi$ in Burn loops}
In this chapter, the kernel of the map $\Phi$ will be determined, for the
case that the loop is of Burn type. The elements of $\ker\Phi$
are of the form $(\lambda, id)$, with $\lambda \in G(L)$; thus $\ker
\Phi$ is isomorphic to a subgroup of $G(L)$, let us denote this
subgroup by $K$. (By Theorem 3.1 of \cite{Funk-Nagy}, even $K\leq
S(N_\lambda)$ holds.)

If $a_1, \ldots, a_k$ are elements of a group, then $[a_1, \ldots, a_k]$
denotes the commutator $a_1^{-1} \cdots a_k^{-1} a_1 \cdots a_k$. Let $L$
be a Burn loop. For $k \geq 2$, we define the following subgroup $H_k$ of
$G(L)$:
\[H_k = \langle [\lambda_{x_1}, \ldots, \lambda_{x_k}] | x_1, \ldots,
x_k \in L, \lambda_{x_1} \cdots \lambda_{x_k} \in S(L)
\rangle.\]
\begin{lm} In any Bol loop, $K = \cup_k H_k$. If the loop is of
Burn type, we have  $\ker \Phi \vartriangleleft G(L)$.
\end{lm}
\pf An element of $\ker \Phi$ is of the form
$(\rho_{x_0} \lambda_{x_0} \cdots \rho_{x_k} \lambda_{x_k},
\lambda_{x_0}^{-1} \cdots \lambda_{x_k}^{-1})$,
where $\lambda_{x_0}^{-1} \cdots \lambda_{x_k}^{-1} = id$, $\lambda_{x_0}
= \lambda_{x_1}^{-1} \cdots \lambda_{x_k}^{-1}$. Thus
\[ x_0 \cdot ( \ldots \cdot (x_{k-2} \cdot x_{k-1} x_k) \ldots ) = 1.\]
The Bol property immediately implies that $\rho_x \lambda_x \rho_y
= \rho_{x y} \lambda_x$ for all $x,y \in L$. Then
\begin{eqnarray*}
\rho_{x_0} \lambda_{x_0} \cdots \rho_{x_k} \lambda_{x_k} & = &
\rho_{x_0 \cdot ( \ldots \cdot (x_{k-2} \cdot x_{k-1} x_k) \ldots )}
\lambda_{x_0} \cdots \lambda_{x_k} \\
& = & \lambda_{x_0} \cdots \lambda_{x_k} \\
& = & \lambda_{x_1}^{-1} \cdots \lambda_{x_k}^{-1} \lambda_{x_0} \cdots
\lambda_{x_k} \\
& = & [\lambda_{x_1}, \ldots, \lambda_{x_k}].
\end{eqnarray*}
By the left inverse property, there exists an $x_0 \in L$ such that
$\lambda_{x_0} \cdots \lambda_{x_k} = id$ if and only if $\lambda_{x_1}
\cdots \lambda_{x_k} \in S(L)$. So we have
\[\ker \Phi = \langle [\lambda_{x_0}, \ldots, \lambda_{x_k}] | x_0, \ldots,
x_k \in L, \lambda_{x_0} \cdots \lambda_{x_k} \in S(L) \rangle =
\bigcup_k H_k.\] 
Since in a Burn loop, the set $S(L)$ is invariant under the conjugation with
elements $\lambda_y$, we have $\ker\Phi \triangleleft G(L)$. \qed
As the square of any element of the Burn loop $L$ is in $N_\lambda$, for all
$n\in N_\lambda$, $x,y\in L$, the commutators $[\lambda_n,\lambda_x]$ and
$[\lambda_x^2, \lambda_y]$ belong to $H_2$. Using this we prove the
following lemma.
\begin{lm} \label{kongr}
Let $\alpha_1, \ldots, \alpha_k \in S(L)$ and $\bar\alpha_i \in
S(N_\lambda)$.
\setcounter{szam}{0}\begin{list}{(\roman{szam})}{\usecounter{szam}}
\item $[\alpha_1, \ldots, \alpha_i, \alpha_{i+1}, \ldots, \alpha_k] \equiv
[\alpha_1, \ldots, \alpha_{i+1}, \alpha_i^{\alpha_{i+1}}, \ldots, \alpha_k]
\pmod{H_2}$;
\item $[\alpha_1, \ldots, (\alpha_i \bar\alpha_i), \ldots, \alpha_k] \equiv
[\alpha_1, \ldots, \bar\alpha_i, \alpha_i, \ldots, \alpha_k] \pmod{H_2}$;
\item $[\alpha_1 \cdots \alpha_k, \bar\alpha_i] \in H_2$;
\item $[\alpha_1, \ldots, \alpha_i, \bar\alpha_i, \ldots, \alpha_k] \equiv
[\alpha_1, \ldots, \alpha_k] \pmod{H_2}$.
\item If the element on the right side of the equivalence (i), (ii) or
(iv) is in $H_k$, then the element on the left side is in $H_k$, as well.
\end{list} \end{lm}
\pf (i) We have $\alpha_1 \cdots \alpha_i \alpha_{i+1} \cdots \alpha_k =
\alpha_1 \cdots \alpha_{i+1} \alpha_i^{\alpha_{i+1}} \cdots \alpha_k.$
On the other hand,
\begin{eqnarray*}
\alpha_1^{-1} \cdots \alpha_i^{-1} \alpha_{i+1}^{-1} \cdots
\alpha_k^{-1} & = &
\alpha_1^{-1} \cdots \alpha_{i+1}^{-1} (\alpha_i^{-1})^{\alpha_{i+1}}
[\alpha_i^{\alpha_{i+1}} (\alpha_i^{-1})^{\alpha_{i+1}^{-1}}]
\alpha_{i+2}^{-1} \cdots \alpha_k^{-1} \\
& = &
\alpha_1^{-1} \cdots \alpha_{i+1}^{-1} (\alpha_i^{-1})^{\alpha_{i+1}}
\cdots \alpha_k^{-1}
[\alpha_i^{\alpha_{i+1}} (\alpha_i^{-1})^{\alpha_{i+1}^{-1}}]^\beta,
\end{eqnarray*}
where $\beta = \alpha_{i+2}^{-1} \cdots \alpha_k^{-1} \in S(L)$. Now,
it is sufficient to show that the expression in the square bracket is
an element of $H_2$:
$\alpha_i^{\alpha_{i+1}} (\alpha_i^{-1})^{\alpha_{i+1}^{-1}} =
[\alpha_{i+1}^2, \alpha_i^{-1}]^{\alpha_{i+1}^{-1}} \in H_2$.

(ii) By some similar calculation one can show that
\[ [\alpha_1, \ldots, (\alpha_i \bar\alpha_i), \ldots, \alpha_k] =
[\alpha_1, \ldots, \bar\alpha_i, \alpha_i, \ldots, \alpha_k]
[\alpha_i, \bar\alpha_i]^{\alpha_{i+1}\cdots \alpha_k},\]
and because of $\bar\alpha_i \in S(N_\lambda)$, the last factor is an
element of $H_2$.

\begin{eqnarray*}
\mbox{ (iii) } [\alpha_1\cdots \alpha_k, \bar\alpha_i] &=&
[\alpha_2\cdots \alpha_k, \bar\alpha_i^{\alpha_1}] [\alpha_1,\bar\alpha_i] \\ 
&\equiv& [\alpha_2\cdots \alpha_k, \bar\alpha_i^{\alpha_1}]
\equiv \cdots \equiv [\alpha_k, \bar\alpha_i^{\alpha_1 \cdots \alpha_k}]
\equiv id \pmod{H_2}.
\end{eqnarray*}
\begin{eqnarray*}
\mbox{ (iv) } [\alpha_1, \ldots, \alpha_i, \bar\alpha_i, \ldots, \alpha_k]
& = & [\alpha_1, \ldots, \alpha_k, \bar\alpha_i^{\alpha_{i+1} \cdots
\alpha_k}]\\
& = & [\alpha_1, \ldots, \alpha_k] [\alpha_1\cdots \alpha_k,
\bar\alpha_i^{\alpha_{i+1} \cdots \alpha_k}] \\
& \stackrel{(iii)}{\equiv} & [\alpha_1, \ldots, \alpha_k] \pmod{H_2}.
\end{eqnarray*}
(v) This follows from $H_2 \triangleleft H_k \triangleleft G(L)$. \qed
\begin{prop} \label{hs-1}
Let $L$ be a Burn loop and $\Phi$ and $H_k$ ($k\geq 2$) be defined as in
the beginning of this section and let $s=|L: N_\lambda|$. Then $\ker \Phi =
H_{s-1}$ if $s \geq 3$, and $\ker \Phi = H_2$ if $s= 1 \mbox{ or } 2$.
\end{prop}
\pf Let $B$ be a set of representatives from the cosets of $N_\lambda$ in $L$
such that $1 \in B$. Then any element of $L$ can be written in a unique way
as the product $n b$, with $n\in N_\lambda$, $b\in B$. Let us choose
elements $x_1, \ldots, x_k$, $x_i = n_i b_i$, from $L$ such that
$\lambda_{x_1} \cdots \lambda_{x_k} \in S(L)$. By \ref{kongr} (ii) and (iv), 
$[\lambda_{x_1}, \ldots, \lambda_{x_k}] \equiv 
	 [\lambda_{b_1}, \ldots, \lambda_{b_k}] \pmod{H_2}$.
Applying \ref{kongr} and $b_i^2 \in N_\lambda$ several times, one gets 
$[\lambda_{x_1}, \ldots, \lambda_{x_k}] \equiv 	 [\lambda_{b_1'},
\ldots, \lambda_{b_m'}] \pmod{H_2}$, where $b_1',\ldots ,b_m'$
are different elements of $B\backslash \{1\}$. Moreover,
$\lambda_{x_1} \cdots \lambda_{x_k} \equiv  \lambda_{b_1'}
\cdots \lambda_{b_m'} \pmod{S(N_\lambda)}$, hence
$[\lambda_{x_1}, \ldots, \lambda_{x_k}] \in H_m$, with $m\leq
|B|-1$. \qed
\begin{cor}
If the loop $L$ is a group, then $\ker \Phi \cong H_2 = L'$.
\end{cor}
\begin{lm} \label{ekvik}
Let the subset $B$ of $L$ be defined as before and let us choose elements
$b_1, b_2, b_3 \in B$ such that $b_3 N_\lambda \cdot (b_2 N_\lambda \cdot
b_3 N_\lambda) = N_\lambda$ holds in the quotient loop $L/N_\lambda$. Then
the followings are equivalent.
\setcounter{szam}{0}\begin{list}{(\roman{szam})}{\usecounter{szam}}
\item $\lambda_{b_1} \lambda_{b_2} \lambda_{b_3} \in S(L).$
\item $\lambda_{b_i} \lambda_{b_j} \lambda_{b_k} \in S(L)$ with $\{i,j,k\}
= \{1,2,3\}$.
\item $\lambda_{b_1} \lambda_{b_2} \in S(L).$
\item $\lambda_{b_i} \lambda_{b_j} \in S(L)$ for all $i,j \in \{1,2,3\}.$
\end{list} \end{lm}
\pf (i) $\Rightarrow$ (iii). From $b_3 N_\lambda \cdot (b_2 N_\lambda \cdot
b_3 N_\lambda) = N_\lambda$ we get $\lambda_{b_1} \lambda_{b_2}
\lambda_{b_3} = \lambda_n$, $n \in N_\lambda$. Hence $\lambda_{b_1}
\lambda_{b_2} = \lambda_{b_3^{-1} n} \in S(L)$.

(iii) $\Rightarrow$ (i). The quotient is a Burn loop, thus $b_3 N_\lambda =
b_2 N_\lambda \cdot b_1 N_\lambda$, $b_2 b_1 = b_3 n$, $\lambda_{b_1}
\lambda_{b_2} = \lambda_n \lambda_{b_3}$, and so $\lambda_{b_1}
\lambda_{b_2} \lambda_{b_3} = \lambda_{b_3^2 n} \in S(L)$.

The equivalence (ii) $\Leftrightarrow$ (iv) can be shown in the same manner.
(iv) $\Rightarrow$ (iii) being trivial, we still have to show (i) $\Rightarrow$
(ii). Supposing (i), we have
\[\lambda_{b_2} \lambda_{b_3} \lambda_{b_1} = \lambda_{b_1}^{-1}
\lambda_{b_1} \lambda_{b_2} \lambda_{b_3} \lambda_{b_1} \in S(L)\]
and
\[S(L) \ni \lambda_{b_3}^{-1} \lambda_{b_2}^{-1} \lambda_{b_1}^{-1}=
\lambda_{b_3} \lambda_{n_3} \lambda_{b_2} \lambda_{n_2} \lambda_{b_1}
\lambda_{n_1} =
\lambda_{b_3} \lambda_{b_2} \lambda_{b_1} \lambda_n,\]
with $n_1, n_2, n_3, n \in N_\lambda$, and so $\lambda_{b_3} \lambda_{b_2}
\lambda_{b_1} \in S(L).$ This is enough to imply (ii). \qed
\begin{prop} \label{kerfi}
If $s = |L:N_\lambda| \leq 7$, then $s\in \{1,2,4\}$ and
\[\ker \Phi = [S(N_\lambda), G(L)] = \langle [\lambda_n, \lambda_x] |
n \in N_\lambda, x \in L \rangle.\]
\end{prop}
\pf The quotient $L/N_\lambda$ is a Bol loop of order $s \leq 7$, and so a
group (cf. \cite{Burn1,Burn2}). In $L$, the square of any element is in $N_\lambda$, since
$L/N_\lambda$ is an elementary abelian 2-group, $s\in \{1,2,4\}$. For $s=1
\mbox{ or } 2$ the statement follows directly from \ref{hs-1}. Let us
suppose that $s=4$. If $b_1 N_\lambda, b_2 N_\lambda, b_3 N_\lambda$ are
different nontrivial elements of $L/N_\lambda$, then $b_3 N_\lambda \cdot
b_2 N_\lambda \cdot b_1 N_\lambda = N_\lambda$. Suppose that $\lambda_{b_1}
\lambda_{b_2}$ or $\lambda_{b_1} \lambda_{b_2} \lambda_{b_3}$ is an elements
of $S(L)$. Then, by \ref{ekvik}, for all $i, j \in \{1,2,3\}$, one has
$\lambda_{b_i} \lambda_{b_j} \in S(L)$. This means that for any $x_i, x_j
\in L$, $x_{i,j} = b_{i,j} n_{i,j}$ with $n_{i,j} \in N_\lambda$,
\[\lambda_{x_i} \lambda_{x_j} = \lambda_{n_i} \lambda_{b_i} \lambda_{n_j}
\lambda_{b_j} = \lambda_{n_j^{T(b_i)} n_i} \lambda_{b_j b_i} \in S(L),\]
thus $L$ is a group, which contradicts $s=4$.

This shows that $\ker \Phi = H_3 = [S(N_\lambda), G(L)]$. \qed
\section{The groups generated by the Bol reflections and the cores of the
loops $B_{4n}$ and $C_{4n}$}
Let us denote by $\sigma_m$ the Bol reflection with axis $x=m$ (see
\cite{Funk-Nagy}), by $N^+$ the collineation group generated by these
reflections and by $N$ the subgroup generated by products of even length of
reflections. Since a Bol reflection interchanges the horizontal and
transversal directions, $N^+$ does not preserve the directions, but the
group $N$ does. 

Clearly, $N$ is a normal subgroup of index 2 of $N^+$ 
and by the geometric properties of Bol reflections, the set $\Sigma =\{
\sigma_x| x \in L\}$ is invariant in $N^+$. Thus, the elements
$\sigma_x\sigma_1$ generate $N$. Using the coordinate system, we get the
form $\sigma_x \sigma_1 = (p_x, \lambda_x)$ for these generators, where 
$p_x=\lambda_x^{-1}\rho_x^{-1}$, see \cite{Ngnote}. 

The following lemma will help us to determine the orbit of the $y$-axe
under $N$. 
\begin{lm} \label{yorbit}
Let $(L,\cdot)$ be a Burn loop and let us define the groups
\[F=\langle p_x | x \in L\rangle, \hskip 1cm
U = \langle \lambda_x^2 | x \in L\rangle.\]
Then, the orbits $1^F$ and $1^U$ coincide. 
\end{lm}
\pf Using that $L$ is left conjugacy closed, we have
\[1^{p_{y_1}\ldots p_{y_k}} = 1^{\lambda_{y_k}^{-1} \ldots
\lambda_{y_1}^{-2} \ldots \lambda_{y_k}^{-1}} = 1^{ \lambda_{y_1'}^{-2}
\ldots \lambda_{y_k'}^{-2}} \in 1^U,\]
which means $1^F \subseteq 1^U$. On the other hand, 
\[1^{p_{y_1}\ldots p_{y_k} \lambda_z^2} = 
1^{\lambda_z \lambda_{y_k'}^{-1} \ldots \lambda_{y_1'}^{-2} \ldots
\lambda_{y_k'}^{-1} \lambda_z}
= 1^{ p_{y_1'} \ldots p_{y_k'} p_z^{-1}} \in 1^F\]
shows that $1^F$ is invariant under $U$. Thus, $1^F=1^U$. \qed
\begin{lm} \label{abelLam}
Let $(L,\cdot)$ be a Burn loop and $U\subseteq G(L)$ be an Abelian group
containing the left translations $\{\lambda_m : m \in N_\lambda\}$. Then
the group $\Phi^{-1}(U)$ of collineations is Abelian, too.
\end{lm}
\pf The action of an arbitrary collineation $(u,v)$ on the set of
transversal lines is $v \lambda_a$, where $a=1^u$, see \cite{BarlStr}. If
$(u,v) \in \Phi^{-1}(U)$, then by Lemma \ref{yorbit} $a \in N_\lambda$,
hence $\lambda_a \in U$ and $v \lambda_a \in U$. And since $U$ is Abelian,
this means that the commutator elements of $\Phi^{-1}(U)$ act trivially on
the set of horizontal and vertical lines, thus on the whole point set. \qed
\begin{thm} Let $(L,\cdot)$ be one of the loops $B_{4n}$ or $C_{4n}$, $n\geq
2$. 
\begin{enumerate}
\item The group $N$ is equal to $\ker\Phi \rtimes \bar{G}$, where $\Phi$
induces an isomorphism from the subgroup $\bar{G}$ to $G(L)$. Let us denote
the respective generators of $\bar{G}$ by $\bar\alpha$, $\bar\beta$ and
$\bar\gamma$, and by $\delta$ the generator of $\ker\Phi$. Then,
$\bar\alpha$ and $\bar\gamma$ act trivially on $\ker\Phi$, and
$\bar\beta \delta \bar\beta =\delta^{-1}$. 
\item The reflection $\sigma_1$ is an automorphism of $N$, which inverts
the generators $(p_x,\lambda_x)$. It always leaves $\bar\alpha$ and
$\bar\beta$ invariant and acts on $\bar\gamma$ and $\delta$ in the
following way.
\[ \sigma_1: \left \{
\begin{array}{l l l} 
\bar\gamma\mapsto \bar\gamma, & \delta \mapsto \bar\alpha^{-4} \delta^{-1} 
& \mbox{if } L=B_{4n}, n\geq 2 \mbox{ or } C_{4n}, n\equiv 0 \pmod{4};
\\
\bar\gamma\mapsto \bar\alpha^n \bar\gamma, & \delta \mapsto \bar\alpha^{-4}
\delta^{-1}
& \mbox{if } L=C_{4n}, n\equiv 2 \pmod{4};
\end{array} \right. \]
\item The group $G_{\mathrm{core}}$ generated by the core is isomorphic to
$N^+/Z(N^+)$ where
\[ Z(N^+) = \left \{ \begin{array}{l @{\mbox{ if }} l} 
\langle \bar\alpha^n, \bar\gamma, \sigma_1 \rangle & L=B_8; \\ 
\langle \bar\alpha^n, \bar\gamma \rangle & L=B_{4n}, n \not\equiv 0
\pmod{4}, n>2; \\ 
\langle \bar\alpha^n, \bar\gamma, \delta^{\frac{n}{4}} \rangle 
& L=B_{4n},  n\equiv 0 \pmod{4}; \\
\langle \bar\gamma \bar\alpha^{\frac{n}{2}}, \delta^{\frac{n}{4}} \rangle
& L=C_{4n}, n\equiv 0 \pmod{4}; \\
\langle \bar\alpha^n \rangle 
& L=C_{4n}, n\equiv 2 \pmod{4}.
\end{array} \right . \]
\end{enumerate}
\end{thm}
\begin{table}
\[ \begin{array}{|c|c|c|c|c|}
\hline
& \begin{array}{c} B_{4n},\\ \mbox{$n$ odd} \end{array}
& \begin{array}{c} B_{4n},\\ \mbox{$n$ even} \end{array}
& \begin{array}{c} C_{4n},\\ n\equiv 2 \pmod{4} \end{array}
& \begin{array}{c} C_{4n},\\ n\equiv 0 \pmod{4} \end{array} \\ 
\hline
\ker \Phi & C_n & C_\frac{n}{2} & C_\frac{n}{2} & C_\frac{n}{2} \\
\hline
|\mbox{($y$-axe)}^N| & n & \frac{n}{2} & n & \frac{n}{2} \\
\hline
\end{array} \]
\caption{The kernel of $\Phi$ and the orbit of the $y$-axe under $N$}
\label{kerphis} \end{table}
\pf If $L$ is either $B_{4n}$, $n\geq 2$ or $C_{4n}$, $n \equiv 0
\pmod{4}$, then by Table \ref{kerphis}, $\ker\Phi$ acts regularly on the
orbit $(\mbox{$y$-axe})^N$. Hence, in these cases, $\bar{G} =
N_{\mbox{$y$-axe}}$ is a good choice. 

Let us suppose $L=C_{4n}$, $n \equiv 2 \pmod{4}$. Let $m$ be
$1^{\alpha^2}$. Then $m$ has order $n$ in $L$, it is a generator of the
cyclic group $N_\lambda$, and the generating element $\delta$ of $\ker \Phi$
can be assumed to be in the form $(\lambda_m^{-2}, id)$. Let $X$ be the set
of vertical lines of equation $x=1$ or 
$x=m^{\frac{n}{2}}$. Let us define the subgroup $\bar{G}$ as the setwise
stabilizer of $X$ in $N$. To the left translation $\lambda_x= \beta\gamma$ the
$N$-generator $(p_x, \lambda_x)$ is associated; since $1^{p_x} = 1^{(\beta
\gamma)^2} = 1^{\alpha^n} =m^{\frac{n}{2}}$, this generator interchanges the
lines in $X$. Therefore $|\bar{G}:N_{\mbox{$y$-axe}}| = 2$ and $|N:\bar{G}|
=n/2$.  Clearly, $\bar{G} \cap \ker\Phi =\{id\}$, and so, $\bar{G}$ is a
transversal to $\ker\Phi$. 

To complete the proof of point 1, we consider the action of $\bar{G}$ on
$\ker\phi$. Applying Lemma \ref{abelLam} to $U=\langle \alpha, \gamma
\rangle$ we see that $\bar\alpha$ and $\bar\gamma$ commute with $\ker\Phi$.
Furthermore, since in each cases of $L$, $\bar\beta \in
N_{\mbox{$y$-axe}}$, hence $\bar\beta = (\beta,\beta) \in N_{(1,1)}$ and
$\delta^{\bar\beta} = \delta^{-1}$. 

\begin{table}
\[ \begin{array}{|l|c|c|}
\hline
(L,\cdot) & \hskip 1cm \lambda_x \hskip 1cm & \hskip 8mm 
(p_x, \lambda_x) \hskip 8mm
\\ \hline
\mbox{(a) } B_{4n},\hskip 2mm C_{4n},\hskip 2mm n \geq 2 & \alpha^{2i} & 
\bar\alpha^{2i} \delta^i 
\\ \cline{2-3}
& \alpha^{2j+1}\beta & \bar\alpha^{2j+1}\bar\beta
\\ \hline
\mbox{(b) } B_{4n},\hskip 2mm n \geq 2 &
\alpha^k\beta\gamma & \bar\alpha^k\bar\beta\bar\gamma
\\ \hline
\mbox{(c) } C_{4n},\hskip 2mm n\equiv 0 \pmod{4} &
\alpha^k\beta\gamma & \bar\alpha^k\bar\beta\bar\gamma \delta^{\frac{n}{4}} 
\\ \hline
\mbox{(d) } C_{4n},\hskip 2mm n\equiv 2 \pmod{4} &
\alpha^k\beta\gamma & \bar\alpha^k\bar\beta\bar\gamma
\\ \hline
\end{array} \]
\caption{Generating elements for $G(L)$ and $N$} 
\label{gens} \end{table}
To determine the action of $\sigma_1$ on the elements $\bar\alpha$,
$\bar\beta$, $\bar\gamma$ and $\delta$, we have to express the generators
$(p_x,\lambda_x)$ of $N$ by these elements. We claim that this is done in
Table \ref{gens}. We therefore use the fact that two collineations $(u,v)$
and $(u',v')$ coincide if $v=v'$ and $1^u=1^{u'}$, see \cite{BarlStr}.
Moreover, if $(u,v)$ is a generator element for $N$, then we have
$1^u=1^{v^{-2}}$. 

Again, the cases $L=B_{4n}$, $n\geq 2$ or $C_{4n}$, $n \equiv 0
\pmod{4}$ are trivial, since then $\bar\alpha$, $\bar\beta$ and
$\bar\gamma$ stabilize the $y$-axe and $\delta$ acts on it in a well known
way. Let us suppose $L=C_{4n}$, $n \equiv 2 \pmod{4}$ and denote the
$N$-generator associated to $\alpha^k \beta\gamma$ by $(u, \alpha^k \beta
\gamma)$. Then one has $1^u = 1^{(\alpha^k \beta \gamma)^2} = 1^{\alpha^n} =
m^{\frac{n}{2}}$, and so, $(u,\alpha^k \beta \gamma) \in \bar{G}$. This gives
$(u,\alpha^k \beta \gamma) = \bar\alpha^k \bar\beta \bar\gamma$. 
The results of Table \ref{gens} and point 2 of the theorem follow. 

The core of the Bol loop $(L,\cdot)$ is the grupoid $(L,+)$ with
$x+y=x\cdot y^{-1}x$. Isomorphic versions of the grupoid can be defined in
the following ways.
\[ \begin{array}{ll} 
(S(L), \oplus), & \lambda_x\oplus \lambda_y =\lambda_x \lambda_y^{-1}
\lambda_x; \\
(\Sigma, \otimes), \hskip 3mm \Sigma=\{\sigma_x:x \in L\}, \hskip 1cm &
\sigma_x\otimes \sigma_y=\sigma_x \sigma_y \sigma_x. 
\end{array} \]
The isomorphism $(L,+) \cong (S(L),\oplus)$ is trivial, and $(S(L), \oplus)
\cong (\Sigma,\otimes)$ can be shown using $\sigma_x\sigma_1 = (p_x,
\lambda_x)$. Hence, the permutation group generated by the core acts on $L$
like $N^+$ acts on $\Sigma$ by conjugation and this action equals to the
action of $N^+$ on the set of vertical lines. And since $\Sigma$ generates
$N^+$, the group $G_{\mathrm{core}}$ generated by the core is isomorphic to
$N^+/Z(N^+)$. Thus, we only have to compute the centres $Z(N^+)$. 

If $L=B_8$, then $\sigma_1$ acts trivially on $N$. In any other cases,
$\sigma_1$ is a non-trivial outer automorphism and we have $Z(N^+) =
C_{Z(N)}(\sigma_1)$, which is very easy to calculate. \qed
\section{Automorphisms of Burn loops of type $B_{4n}$ and $C_{4n}$}
Let $(L,\cdot)$ be a loop and let $u$ denote an automorphism of
$L$. Then, by conjugation, $u$ induces an automorphism of the
group $G(L)$. Moreover $u$ leaves the section $S(L)$ and the
stabilizer $G(L)_1$ invariant.
Conversely, let $u$ be an automorphism of $G(L)$, normalizing
the subgroup $G(L)_1$ and the set $S(L)$. Then $u$ induces a
permutation on the cosets of $G(L)_1$, hence on $L$. The induced
permuation will fix 1 and normalize $S(L)$, thus $u^{-1}\lambda_x u
= \lambda_y$ for all $x \in L$.  Applying this to 1, one gets
$y=x^u$, hence $\lambda_x^u = \lambda_{x^u}$ for all $x \in L$.
This means $u \in \aut (L)$.

In the case of the given loops the stabilizer of 1 consists of $\{ id,
\beta\}$. First we calculate its normalizer in the automorphism
groups of the left translation groups, that is, the groups
$C_{\aut (G)}(\beta)$, where $G$ is $G_{8n}$ or $H_{8n}$. 

\begin{lm} \label{autg8nodd}
Let $G$ denote the group $G_{8n}$, $n$ odd. Then
$C_{\aut (G)}(\beta) \cong Z_{2n}^\ast \times S_3$, 
and the elements of $C_{\aut (G)}(\beta)$ normalize $S(B_{4n})$.
\end{lm}
\pf Let us define the subgroups $A=\langle \alpha^2 \rangle$ and
$B=\langle \alpha^n, \beta, \gamma \rangle$ of $G$. As $|A|=n$ is
odd, $A$ is a characteristic subgroup of $G=A\times B$. Moreover, $B=Z(G)
\langle \beta \rangle$ is invariant in $C_{\aut (G)}(\beta)$, as
well. Hence, $C_{\aut (G)}(\beta)= \aut (A) \times C_{\aut (B)}(\beta)
\cong Z_n^\ast \times S_3$. 

On the other hand, $S(L)=A\{id, \alpha^n \beta, \beta\gamma,
\alpha^n \beta\gamma\}$. Since the set 
\[\{id, \alpha^n \beta, \beta\gamma, \alpha^n \beta\gamma\}\]
is invariant under $C_{\aut (B)}(\beta)$, the statement follows. \qed
\begin{lm} \label{autg8neven}
Let $G$ denote the group $G_{8n}$, $n$ even. Then 
$C_{\aut (G)}(\beta) \cong Z_n^\ast \times D_8$,
and the elements of $C_{\aut (G)}(\beta)$ normalize
$S(B_{4n})$. 
\end{lm}
\pf It is enough to consider the possible images of $\alpha$ and
$\gamma$, let us write them as $\hat\alpha = \alpha^i \gamma^k
\beta^j$ and $\hat\gamma = \alpha^p \gamma^q \beta^s$,
respectively. Clearly, $\hat\beta=\beta$.

If $j=1$ then $\hat\alpha^2=id$, which is impossible. The order
of $\hat\alpha$ must be $2n$, thus $i \in Z_{2n}^\ast$. The
elements $\hat\alpha$ and $\hat\gamma$ must commute, $s$ cannot
be 1. Also the elements $\hat\beta$ and $\hat\gamma$ commute, we
must have $p=ln$ with $l\in Z_2$. 

Let us now suppose that $q=0$. Then $l=0$ implies
$\hat\gamma=id$ and $k=0$ implies $\gamma \not \in \langle
\hat\alpha, \hat\beta, \hat\gamma \rangle$, hence we have
$l=k=1$. This means $\hat\alpha^n = \alpha^{ni}=\alpha^n = \hat
\gamma$, a contradiction. 

Let us denote by $u(i,k,l)$ the automorphism induced by 
\[\alpha\mapsto \alpha^i \gamma^k, \hskip 1cm
\beta \mapsto \beta, \hskip 1cm
\gamma\mapsto \alpha^{ln} \gamma,\]
with $i \in Z_{2n}^\ast$, $k,l \in Z_2$. It is easy to check
that this is really an element of $C_{\aut (G)}(\beta)$.
Moreover, 
\[u(i,j,k) u(i',j',k') = u(ii'+lk'n, k+k',l+l'),\]
where one calculates modulo $2n$ in the first and modulo 2 in
the second and third position. 

Let us decompose $Z_{2n}^\ast$ into $Z_n^\ast \times Z_2$ by
$i=i_0+i_1n$, $i_0\in Z_n^\ast$, $i_1 \in Z_2$. Then the group
$C_{\aut (G)}(\beta)$ decomposes into the direct factors
\[\{u(i_0,0,0) : i_0 \in Z_n^\ast \} \mbox{ and } 
\{u(i_1 n, k, l) : i_1,k,l \in Z_2\}.\]
An easy calculation is to show that the second factor is
isomorphic to the dihedral group $D_8$ of 8 elements. 

Since we gave explicitely the elements of $C_{\aut (G)}(\beta)$,
it can be checked directly that they leave $S(L)$ invariant. \qed
\begin{lm} \label{auth8n}
Let $G$ denote the group $H_{8n}$, $n>2$ even. Then 
$C_{\aut (G)}(\beta) \cong Z_{2n}^\ast \times Z_2$,
and the elements of $C_{\aut (G)}(\beta)$ normalize
$S(C_{4n})$. 
\end{lm}
\pf As in the preceding proof, we consider the images $\hat
\alpha=\alpha^i\gamma^k \beta^j$, $\hat\gamma= \alpha^p \gamma^q
\beta^s$ of $\alpha$ and $\gamma$. 

If $j=1$, then $\hat\alpha^2=\alpha^i\gamma^k \beta
\alpha^i\gamma^k \beta = (\gamma^k \beta)^2 = \alpha^{kn}$,
$\hat\alpha^4=id$, which is not possible because of $n>2$. 
If $k=1$, then $(\hat\alpha \hat\beta)^2=(\gamma\beta)^2 =
\alpha^n \neq id$, hence $k=0$ and $\hat\alpha=\alpha^i$, with
$i \in Z_{2n}^\ast$. 

As before, $\hat\alpha \hat\gamma = \hat\gamma \hat\alpha$
implies $s=0$ and $\gamma \in \langle \hat\alpha, \hat\beta,
\hat\gamma \rangle$ implies $q \neq 0$. Finally, $p\in \{0,n\}$,
since $\hat\gamma = (\alpha^p \gamma)^2 = \alpha^{2p} =id$.

Thus, any element of $C_{\aut (G)}(\beta)$ is induced by
\[\alpha\mapsto \alpha^i, \hskip 1cm
\beta \mapsto \beta, \hskip 1cm
\gamma\mapsto \alpha^{ln} \gamma,\]
and it leaves $S(L)$ invariant. \qed
\begin{thm} \label{loopaut}
Let $(L,\cdot)$ be one of the loops $B_{4n}$ or $C_{4n}$ defined
at the beginning of this section. Then
\[\aut (L) \cong \left \{
\begin{array}{ll}
Z_n^\ast \times S_3 & \mbox{if $L=B_{4n}$, $n$ odd} \\
Z_n^\ast \times D_8 & \mbox{if $L=B_{4n}$, $n$ even} \\
Z_{2n}^\ast \times Z_2 & \mbox{if $L=C_{4n}$, $n>2$, $n$ even} \\
D_8 & \mbox{if $L=C_8$}
\end{array} \right. \]
Moreover, in any of these loops, each left pseudo-automorphism
is an automorphism. 
\end{thm}
\pf The case $L=C_8$ is handled in \cite{Ngnote}, the others in Lemmas
\ref{autg8nodd}, \ref{autg8neven} and \ref{auth8n}. We only have
to prove the second statement. Therefore, let us suppose that
$u$ is a left pseudo-automorphism of $L$ with companion $c$,
that is, for all $x,y \in L$,
\[(c\cdot x^u)\cdot y^u = c\cdot (xy)^u.\]
This can be expressed by $u \lambda{c x^u} = \lambda_x u
\lambda_c$, which implies $S(L)^u = S(L) \lambda_c^{-1}$. 

The following results are to find in \cite{Burn1,Burn2}.
If $L=B_{4n}$, then the principal isotopes of $L$ have
the four representation $S(L)$, $\alpha \beta S(L)$,
$\alpha\beta \gamma S(L)$, and $\beta \gamma S(L)$. If $n$ is
even, then these sections contain $3n+1$, $n+3$, $n+3$ and $n+1$
elements of order 2.
If $n$ is odd, $S(L)$ contains $3n$ elements of order 2 and the
others contain $n+2$ elements of order 2, $n>2$.  
That means that $c$ is a left companion element of $L$
if and only if $S(L) \lambda_c = S(L)$, it is, $c \in
N_\lambda$ and $u$ is an automorphism. 

Let now $L$ be equal to $C_{4n}$. Again the principal isotopes are 
$S(C_{4n})$, $\alpha \beta S(C_{4n})$, $\alpha \beta
\gamma S(C_{4n})$, and $\beta \gamma S(C_{4n})$, they contain
$n+1$, $n+3$, $3$ and $1$ involutions, respectively. If $n>2$,
then one sees with the above argument that $c \in
N_\lambda$ and $u$ is an automorphism. \qed
\section{Collineation groups of the given 3-nets}
In this chapter, we determine the full collineation group $\Gamma$ of the
3-nets belonging to $B_{4n}$, $n\geq 3$, and $C_{4n}$, $n\geq 4$, $n$ even.
The cases $B_8$ and $C_8$ are completely descripted in \cite{Ngnote}.

Denote by $P$ the orbit $(1,1)^\Gamma$ of the origin under
$\Gamma$. As we know by Corollary 2.8 of \cite{Ngnote}, for any
Burn loop, $P$ is a union of vertical lines and its intersection
with the $x$-axe constitute of the points belonging to the left
companion elements. In our cases, these are the elements of
$N_\lambda$, see Theorem \ref{loopaut}. Hence $|P|= 4n^2$.

Let $\Lambda_0$ be the subgroup $\langle \alpha, \gamma \rangle$
of $G(L)$. The centralizer element $\alpha^i \beta \gamma^j \not \in \Lambda_0$
in $\Lambda_0$ has order 4, that is, any Abelian subgroup not contained in
$\Lambda_0$ has order at most 8. This means that if $n>2$ then $\Lambda_0$ is
the only Abelian subgroup of index 2 in $G(L)$, it must therefore be
characteristic in $G(L)$. 

Now, we define the following subgroups of $\Gamma$. 
\[ \begin{array}{ll}
T = \{ (\lambda_m,id) : m\in N_\lambda\},
\hskip 1cm & \Lambda = \Phi^{-1}(\Lambda_0), \\
A = \{ (\sigma,\sigma) : \sigma \in \aut (L)\}, \hskip 1cm &
M=T\Lambda.
\end{array} \]
\begin{lm}\label{Mlemma}
The subgroup $M$ is an Abelian normal subgroup of $\Gamma$. Moreover, it is
isomorphic to the direct product $N_\lambda \times \Lambda_0$ and acts
regularly on the orbit $P$ of the origin.
\end{lm}
\pf First we show that $M$ is Abelian. By Lemma \ref{abelLam}, one sees
that the permutation action of the elements of $\Lambda$ are all in $\langle
\alpha, \gamma \rangle$; the same can be said about the elements of $T$. 
These actions commute, and so, all the elements must commute. 

Clearly, $T$ is normal in $\Gamma$. The subgroup $\Lambda$ is invariant in
$\Gamma$ as well, for it is the homomorphic preimage of a characteristic
subgroup.  

Suppose that $(u,v)$ is an element of $M_{(1,1)}$. Then $v=id$, since
$v=\beta$ is not possible. This implies $u=\lambda_m$, $m\in N_\lambda$;
from which $u=id$ follows. Furthermore, on the one hand, by $\Lambda\cap T =
\ker \Phi$, we have $M_{\mbox{$y$-axe}} \cong M/T  \cong \Lambda_0$. On the
other hand, $T \subset M_{\mbox{$x$-axe}}$ acts transitively  on $P \cap
\mbox{$y$-axe}$. This means that $M$ acts transitively on
$P$, thus, regularly. Finally, $M = T \times M_{\mbox{$y$-axe}} \cong
N_\lambda \times \Lambda_0$. \qed
\begin{thm}
Let $\Gamma$ be the full collineation group of the 3-net, coordinatized by
the loop $L$, with $L=B_{4n}$ or $C_{4n}$, $n >2$.
Then, $\Gamma$ can be written as the semidirect product $M \rtimes
\aut (L)$, where $M$ is defined as above and the action of $\aut (L)$ on $M$ 
is defined by $(u,v)^{\sigma} = (u^\sigma, v^\sigma)$. 
\end{thm}
\pf Obviously, $A$ is isomorphic to $\aut (L)$. By Theorem 10.1 of
\cite{BarlStr}, $A$ is equal to the stabilizer
$\Gamma_{(1,1)}$ of the origin $(1,1)$ in $\Gamma$. By Lemma 
\ref{Mlemma}, $M$ is a normal subgroup of $\Gamma$, acting
regularly on the orbit $P=(1,1)^\Gamma$. Then, $\Gamma$ can be
written as the semidirect product $M \rtimes A \cong M \rtimes \aut (L)$. 
\qed

{\it Remark.} Note that there is an interesting analogy with the case of
group 3-nets: then one has $\Gamma \cong (G\times G) \rtimes \aut (G)$ (cf.
\cite{BarlStr}, Theorem 10.1). 
\bibliographystyle{plain}
\bibliography{looplit}

\end{document}